\documentclass{amsart}
\usepackage{graphicx}
\usepackage{amssymb}
\newtheorem{thm}{Theorem}[section]

\newtheorem{lem}[thm]{Lemma}

\newtheorem{alg}[thm]{Algorithm}
\newtheorem{exmp}[thm]{Example}

\theoremstyle{definition}
\newtheorem{defn}[thm]{Definition}
\theoremstyle{remark}

\numberwithin{equation}{section}

\begin{document}

\title[Converting  Subalgebra Bases with the Sagbi Walk ]{Converting  Subalgebra Bases with the Sagbi Walk }%
\author{Junaid Alam Khan$^*$}%
\address{$*$Abdus Salam School of Mathematical Sciences, GCU, Lahore Pakistan}%
\email{junaidalamkhan$@$gmail.com}%

\thanks{This research was partially supported by Higher Education Commission, Pakistan  }


\begin{abstract}
We present an algorithm which converts a given Sagbi basis of a polynomial $K$-subalgebra $\mathcal{A}$ to a Sagbi basis of $\mathcal{A}$ in a polynomial ring with respect to another term ordering, under the assumption that subalgebra $\mathcal{A}$ admits a finite Sagbi basis with respect to all term ordering. The Sagbi walk method converts a Sagbi basis by partitioning the computations following a path in the Sagbi Fan. The algorithms have been implemented as a library for the computer algebra system SINGULAR \cite{GPS1}.
\end{abstract}
\maketitle

\section{Introduction and Preliminaries}
Let $K$ be a field and  $K[x_1,\ldots,x_n]$ the polynomial the ring over the field $K$ in $n$ variables.
Let $\mathcal{A}$ be a $K$-subalgebra  of $K[x_1,\ldots,x_n]$. In this paper  we only consider subalgebras
which admits a finite Sagbi basis (see \cite{RS1}) with respect to all global monomial orderings.
The objective of this paper is the presentation of a procedure for converting a Sagbi basis of subalgebra of a subalgebra
$\mathcal{A}$ to a Sagbi basis of $\mathcal{A}$ with respect to another term ordering. We called this procedure Sagbi walk.
This procedure is the Subalgebra analogue of the Gr\"{o}bner walk (see \cite{CKM1}), which used to convert Gr\"{o}bner bases of a Polynomial ideal from one ordering to another ordering.
\\We use the notations from \cite{CLO1} and \cite{GP1} and repeat them for the convenience of reader.
\begin{defn}
A monomial ordering is a total ordering $>$ on the set of monomials $Mon_n=\{x^\alpha \,|\, \alpha \in \mathbb{N}^n\}$ in
n variables satisfying
$$x^\alpha > x^\beta \Longrightarrow x^\gamma x^\alpha > x^\gamma x^\beta$$ for all $\alpha, \beta, \gamma  \in \mathbb{N}^n$. We also say $>$ is a monomial ordering on $K[x_1,\ldots, x_n]$  meaning that $>$ is a monomial ordering on $Mon_n$.
\end{defn}
In this paper we consider only global ordering, i.e $x^{\alpha}> 1$ for all $\alpha\neq 0$.

\begin{defn}
Let $>$ be a fixed monomial ordering. Write $f\in K[x_1,\ldots,x_n]$, $f\neq 0$, in a unique way as a sum of non-zero terms $$f=a_\alpha x^\alpha +a_\beta x^\beta +\ldots +a_\gamma x^\gamma , \, \,\,\,\,\,\,\,\,\, x^\alpha > x^\beta > \ldots  > x^\gamma ,$$ and $a_\alpha ,a_\beta ,\ldots a_\gamma \in K$. We define:
\begin{itemize}
\item[1.]$LM_{>}(f):=x^\alpha$, the leading  monomial of $f$,
\item[2.]$LE_{>}(f):=\alpha$, the leading  exponent of $f$,
\item[3.]$LT_{>}(f):=a_\alpha x^\alpha$, the leading  term of  $f$,
\item[4.]$LC_{>}(f):=a_\alpha$, the leading  monomial of $f$,
\item[7.]$support(f):=\{x^{\alpha},x^{\beta},\ldots,x^{\gamma}\}$, the set of all monomials of $f$ with non-zero coefficient.
\end{itemize}
\end{defn}

\begin{defn}
Let $G$ be any subset of $K[x_1,\ldots,x_n]$ and $>$ be a fixed monomial ordering.\\ \\
We denote by $L_>(G)$, the set of leading monomials of $G$,
$L_{>}(G):=\{\,LM_>(f)\,|\,f\in G \backslash\{0\}\,\}$
\\ \\We define $in_>(G)$ to be a $K$-subalgebra generated by $LT_>(G)$.
$$in_>(G):=K[L_>(G)] $$
\end{defn}
\begin{defn}
Let $>$ is a fixed monomial ordering. A subset $S\subset \mathcal{A}$ is called Sagbi basis of $\mathcal{A}$ with respect to $>$ if
$$in_{>}(\mathcal{A})= K[L_{>}(S)]$$

\end{defn}

\begin{defn}
Let $G$ be a subset of $K[x_1,\ldots,x_n]$ and $>$ is a fixed monomial order.\\

1) $G$ is called interreduced if $0\notin G$ and for any polynomial $f\in G$, $LM_{>}(f)\notin K[LM_{>}(G\backslash \{0\})]$.\\

2) $f\in K[x_1,\ldots,x_n]$ is called  reduced with respect to $G$ if no term of $f$ is contained in $K[LM_{>}(G)]$. \\

3) $G$ is called reduced if $G$ is interreduced, and for any $f\in G$, $LC_{>}(f)=1$ and $tail(f)$ is reduced with respect to $G$.

\end{defn}

\begin{defn}
Given a vector $\mathbf{w}=(w_1,\dots,w_n)\in \mathbb{R}^{n}$, we define the $\mathbf{w}$-$degree$ (weighted degree)  of $x^{\alpha}$ by
$$deg_{\mathbf{w}}(x^{\alpha})=\mathbf{w}.\alpha=w_1\alpha_1+\dots,w_n\alpha_n\,\,,     $$ that is, the variable $x_{i}$ has degree $w_{i}$.\\ \\
 For a polynomial $f=\sum_{\alpha}a_{\alpha}x^{\alpha}$, we define the $w$-$ degree$ by
$$ deg_{w}(f)=max\{\,deg_{w}(x^{\alpha})\,|\,a_{\alpha}\neq 0\,\}$$
The initial forms of $f$ with respect $\mathbf{w}$, denoted by $In_{w}(f)$ is the sum of all those terms in $f$ with maximum of $\mathbf{w}$-$ degree$. Let $deg_{w}(f)=d$ then,
$$In_{w}(f)=\sum_{deg_{w}(x^{\alpha})=d}a_{\alpha}x^{\alpha} $$
Further more we define $deg_{w}(0)=-1$ and $In_{w}(0)=0$.\\ \\
A polynomial $f$ is called $\mathbf{w}$-$homogenous$ if $f=In_{w}(f)$. Note for $w=(1,\dots,1)$ this notion coincide
with standard homogeneity.

\end{defn}
\begin{defn}
Let $G$ be any subset of $K[x_1,\ldots,x_n]$ and $\mathbf{w}$ be a fixed vector.\\ \\
We denote $In_{w}(G)$ as
$$In_w(G):=\{\,In_w(f)\,|\,f\in G \,\}$$
We define $in_w(G)$ to be a $K$-subalgebra generated by $In_w(G)$.
$$in_w(G):=K[In_w(G)] $$
\end{defn}

\begin{defn}
Let $M\in GL(n,\mathbb{R})$. We can use $M$  to obtain a monomial ordering by setting
   $$x^{\alpha}>_{M}x^{\beta}\,:\Longleftrightarrow\,M{\alpha}>M{\beta},$$
where $>$ on the right-hand side is the lexicographical ordering on $\mathbb{R}^{n}$.

\end{defn}

\begin{lem}
(c.f. \cite{GPS1}, page 18 )Any monomial ordering can be defined as $>_{M}$ by a matrix $M\in GL(n,\mathbb{R})$.

\end{lem}

 \begin{defn}
A monomial ordering $>$ on $\{\,x^{\alpha}\,|\,\alpha\in {\mathbb{N}}^{n}\,\}$ is called weighted degree ordering if there
exist a vector $w=(w_1,\ldots,w_n)$ positive integer such that
$$deg_{w}(x^{\alpha})> deg_{w}(x^{\beta})\Rightarrow x^{\alpha}> x^{\beta} $$
Consider $A\in GL(n,\mathbb{R})$. A matrix $>_{A}$ is a global ordering if and only if the first non-zero entry in each column is positive.
It is a weighted degree ordering if and only if all entries in the first row are positive.
\end{defn}

\begin{lem}
(c.f. \cite{GPS1}, page 15 )Let $>$ be a monomial ordering  and $N\subset Mon(x_1,\ldots,x_n)$ a finite set.
Then their exist some $w=(w_1,\ldots,w_n)\in {\mathbb{Z}}^n$ such that $x^{\alpha}> x^{\beta}$ if and only if
$deg_{w}(x^{\alpha})> deg_{w}(x^{\beta})$ for all $x^{\alpha},x^{\beta}\in N $. Moreever, $w$ can be chosen such that
$w_{i}> 0$ if $x_{i}> 1$ and $w_i< 0$ if $x_{i}< 1$.
\end{lem}
\begin{lem}
For $f,g \in K[x_1,\ldots,x_n]$,$\mathbf{w}\in (\mathbb{R}^{n})^{+}$. We have $In_{w}(fg)=In_{w}(f)In_{w}(g)$.
\end{lem}
\begin{defn}
Let $G$ be a subset  of $K[x_1,\ldots,x_n]$
\\
\\ A $G$-monomial is a finite power product of the form $ G^\alpha=g_1^{\alpha_{1}}\ldots g_m^{\alpha_{m}}$ where $g_i\in G$ for $i=1,\ldots,m,$ and $\alpha=(\alpha_1,\ldots,\alpha_m)\in \mathbb{N}^m.$
 \\
\\ The set of all $G$-monomial is denoted by:$$Mon_G=\{G^\alpha| \ \alpha\in \mathbb{N}^m,\,m\in \mathbb{N} \}$$

\end{defn}
The Sagbi Walk breaks the conversion problem into several steps between adjacent Sagbi bases following a path in Sagbi fan
which is the analogue of the Gr\"{o}bner Fan (see \cite{MR1}). Since two term orderings leading to adjacent Sagbi basis can be viewed as
refinement of a common partial ordering, these transformation can be computed working just with the initial form with
respect to this partial ordering.
\section{Sagbi cone}
 Consider subalgebra $\mathcal{A}$  which admits a finite Sagbi basis with respect all monomial ordering.
\begin{lem}
Let $F=\{f_1,\ldots,f_s\}\subset \mathcal{A}$ and $>_1$ and $>_2$ are two monomial order such that $LT_{>_1}(f_i)=LT_{>_2}(f_i)$
for all $i$. If $F$ were is a Sagbi basis of $\mathcal{A}$ with respect to $>_1$ then $F$ is a Sagbi basis of $\mathcal{A}$ with respect to $>_2$.
\end{lem}
\begin{proof}
let $f\in \mathcal{A}$ be arbitrary. Reducing $f$ by $F$ using $>_1$, we obtain
  $$f=a_{1}F^{\alpha_1}+ \dots + a_{s}F^{\alpha_s}+h $$
  where $h$ is Sagbi normal form of $h$ with respect to $F$ using $>_2$, therefore either $h$  is zero or for all
  $ x^{\beta}\in support(h),x^{\beta}\neq LM_{>_2}(F^{\alpha})$ for all F-monomial $F^{\alpha}$. However, we have
   $LT_{>}(f_i)=LT_{>_1}(f_i)$ for all $i$.
  Since $h=f-a_{1}F^{\alpha_1}+ \dots + a_{s}F^{\alpha_s}\in \mathcal{A}$, and $F$ is assumed to be Sagbi basis for $\mathcal{A}$ with respect
   to $>_1$, this implies that $h=0$. In above  we used Sagbi normal form algorithm, therefore
   $LT_{>_2}(f)=LT_{>_2}(a_{i}F^{\alpha_i})$ for some $i$, this shows that $F$ is also Sagbi basis for $\mathcal{A}$ with respect to $>_2$.
\end{proof}

\begin{thm}
Given a subalgebra $\mathcal{A}\subset K[x_1,\ldots,x_n]$ and monomial ordering $>$, there exists a weight
$\mathbf{w}$ such that $in_{>}(\mathcal{A})=in_{w}(\mathcal{A})$.
\end{thm}

\begin{proof}
Let $G=\{g_1,g_2,\ldots,g_m\}$ be a  Sagbi basis with respect to $>$.
Consider all the pairs $(LM_>(g_i),u) $ where $u\in support(g_i)$ and $u\neq LM_>(g_i)$. There are finitely many pairs. Hence by Lemma $10$ there exist a weight $\mathbf{w}$ such that
$deg_{w}(LM_>(g_i))>deg_{w}(u)$ for all $u\in support(g_i)$ with $u\neq LT_>(g_i)$ an for all $i$. This implies specially
$LT_>(g_i)=LT_{>}(g_i)=In_{w}(g_i)$.

\end{proof}

Its is natural to ask following question about the collection of all Sagbi bases of a fixed subalgebra $\mathcal{A}$.\\
$\bullet$ Is the collection of possible Sagbi bases of $\mathcal{A}$ finite or infinite$?$\\
$\bullet$ When do different monomial ordering yields the same monic (reduced) Sagbi basis for $\mathcal{A}$?\\

Given a subalgebra $\mathcal{A}\subset K[x_1,\ldots,x_n]$, we define the set $$IN(\mathcal{A})=\{\,in_>(\mathcal{A})\,|\,>\, \hbox{ a monomial ordering}\,\}$$\\
Now we answer the first question

 \begin{thm}
 For a subalgebra $\mathcal{A}\subset K[x_1,\ldots,x_n]$ the set $IN(\mathcal{A})$ is finite.

 \end{thm}

\begin{proof}
Aiming for a contradiction, suppose that $IN(\mathcal{A})$ is infinite set. For each $N\in IN(\mathcal{A})$, let $>_n$ be any one particular
monomial order such that $N=in_{>_n}(\mathcal{A})$. Let $\Gamma$ be collection of monomial orders $\{\,>_n\,|\,N\in IN(\mathcal{A})\,\}$. Our assumption implies that $\Gamma$ is infinite.

As $\mathcal{A}$ is finitely generated subalgebra we have $\mathcal{A}=K[f_1,\ldots,f_s]$ for polynomials
$f_i\in  K[x_1,\ldots,x_n]$. Since each $f_i$ contains only a finite number of terms, by a pigeonhole principle argument, there exist a infinite subset $\Gamma_1\subset \Gamma$ such that the leading terms $LT_{>}(f_i)$ agree for all $>$ in $\Gamma_1$ and all $i,\, i\leq i\leq s$. We write $N_1=K[LT_>(f_1),\ldots,LT_>(f_s)]\,$
( taking any monomial order $>$ in $\Gamma_1$).

If $F=\{\,f_1,\dots,f_s\}$ were a Sagbi basis for $\mathcal{A}$ with respect to some $>_1$ in $\Gamma_1$,
then by lemma 13 $F$
 would be a Sagbi basis for $\mathcal{A}$ with respect to every $>$ in $\Gamma_1$.

  However, this cannot be the case since original set of monomial orders $\Gamma\supset \Gamma_1$ was chosen so that the
  $in_>(\mathcal{A})$ for $>$ in $\Gamma$ were all distinct. Hence, given any $>_1$ in $\Gamma$, there must be some $f_{s+1}\in \mathcal{A}$
  such that $LT_{>_1}(f_{s+1})\notin K[LT_{>_1}(f_1),\dots,LT_{>_1}(f_s)]=N_1$. Replacing $f_{s+1}$ by its Sagbi normal form
  on reducing by $f_1,\dots,f_s$, we may assume that for all $ x^{\beta}\in support(f_{s+1}),x^{\beta}\neq LM(F^{\alpha})$
  for all F-monomial $F^{\alpha}$ where $F=\{\,f_1,\dots,f_s \,\}$.

  Now apply pigeonhole principle again to find an infinite subset $\Gamma_2\subset \Gamma_1$ such that the leading terms of
   $f_1,\dots,f_{s+1}$ are the same for all $>$ in $\Gamma_2$. Let $N_2=K[LT_{>_1}(f_1),\dots,LT_{>_1}(f_{s+1})]$ for all
   $>$ in $\Gamma_2$, and note that $N_2\subset N_2$. the argument in the preceding paragraph shows that
   $\{f_1,\dots,f_{s+1}\}$ cannot be a Sagbi basis with respect to any monomial order in
   $\Gamma_1$, so fixing $>_2\in \Gamma_2$, we find an $f_{s+2}\in \mathcal{A}$ such that for all
   $ x^{\beta}\in support(f_{s+2}),x^{\beta}\neq LM(F^{\alpha})$ for all F-monomial $F^{\alpha}$
   where $F=\{\,f_1,\dots,f_{s+1} \,\}$.

  Continuing in this way, we produce a descending chain of infinite subsets
  $\Gamma\supset \Gamma_1\supset \Gamma_2\supset \Gamma_3\dots,$ and an infinite strictly ascending chain
  $N_1\subset N_2 \subset N_3\subset\dots$. This contradicts the  condition that $\mathcal{A}$ admits a finite Sagbi basis. so the proof is complete.
\end{proof}
Let $S=\{\,g_1,\dots,g_t\,\}$ be one of the Sagbi bases of $\mathcal{A}$ with respect to a monomial ordering $>$ such that $LT(g_i)=x^{\alpha(i)}$, and $N=K[x^{\alpha(1)},\dots,x^{\alpha(t)}]$ the corresponding element of $In(\mathcal{A})$. Our next goal is to understand  the set of monomial ordering for which $S$ is the corresponding Sagbi basis of $\mathcal{A}$. This will answer the second question posed at the start of this section. we write
 $$ g_i=x^{\alpha(i)}+\sum\limits_{\beta}c_{i,\beta}x^{\beta}\,\,,$$
 where $x^{\alpha(i)}>x^{\beta}$ whenever $c_{i,\beta}\neq 0$. by the above discussion, each such ordering $>$ comes from a matrix $M$, so in particular to find the leading terms, we  compare monomials first according to first row $\mathbf{w}$ of the matrix.

 If $\alpha(i).\mathbf{w }>{\beta}.\mathbf{w}$ for all $\beta$ with $c_{i,\beta}\neq 0$, the single weight vector $\mathbf{w}$ selects the correct leading term in $g_i$ as the term of highest of weight. As we know, however, we may have a tie in the first computation using other rows of $M$. There may be also some other weight vector $\mathbf{w}$
 which qualify to become first row of $M$.  This suggest that we should define certain set of weight vectors.
 \begin{defn}
Let $S$ be a  Sagbi basis of $\mathcal{A}$ with respect to a monomial ordering $>$ as in above discussion.
We define the $C_{S,>}$ as
 $$ C_{S,>}=\{\,\mathbf{w}\in {(\mathbb{R}^{n})}^{+}\,:\,\alpha(i).\mathbf{w}\geq \beta.\mathbf{w}\,\,\,\hbox{whenever}\,\,\,c_{i,\beta}\neq 0\,\}$$

$$ \hspace{1.1 cm}=\{\,\mathbf{w}\in {(\mathbb{R}^{n})}^{+}\,:\,(\alpha(i)- \beta).\mathbf{w}\geq 0\,\,\,\hbox{whenever}\,\,\,c_{i,\beta}\neq 0\,\}$$
$C_{S,>}$  is called \textbf{Sagbi cone}\footnote{  If $S$ is a Grobner basis of some ideal $I$ with
respect to $>$, then $C_{S,>}$ is called Gr\"{o}bner Cone (see \cite{CLO1}). It is also defined as the  closure in $\mathbb{Q}^{n}$ of
$\{w\in (\mathbb{Q}^{n})^{+}\,|\,\langle LT_{>}(I)\rangle=\langle In_{w}(I)\rangle\,\}$(see \cite{CKM1}.
Similarly for
subalgebra $\mathcal{A}$ we can define the Sagbi cone
as the  closure in $\mathbb{Q}^{n}$ of  $\{w\in (\mathbb{Q}^{n})^{+}\,|\, in_{>}(\mathcal{A})= in_{w}(\mathcal{A})\,\}$.
By theorem $2.14$ its interior is nonempty so it is well defined. }
it is a closed,convex polyhedral cone. \footnote{ Sagbi cones have same geometrical structure as Gr\"{o}bner cones.
For Gr\"{o}bner cones properties see \cite{CLO1}.   }.
By lemma $8$, its interior\footnote{ Its interior is defined as int$(C_{S,>})=\{\,\mathbf{w}\in {(\mathbb{R}^{n})}^{+}\,:\,\alpha(i).\mathbf{w}> \beta.\mathbf{w}\,\,\,\hbox{whenever}\,\,\,c_{i,\beta}\neq 0\,\} $  } is non-empty. Note that for $>_M$, a matrix order such that first row of M lies in the int($C_{S,>}$).
 $S$ is also Sagbi basis with respect to $>_{M}$. It gives an answer to the second question.

The collection of all the cones $C_{G,>}$ and faces\footnote{A face of a cone $\sigma$ is $\sigma \cap \{l=0\}$,
where $l=0$ is a non-trivial linear equation such that $l\geq 0$ on $\sigma$.  } as $ G$ ranges over all  Sagbi basis
 with respect to $>$ of
$\mathcal{A}$ is called \textbf{Sagbi Fan}. This collection is finite follows from theorem acoording to theorem 2.6 .
 \end{defn}
\section{Sagbi Walk}
Let assume that  $S$ is a Sagbi basis of $\mathcal{A}$ with respect to some monomial ordering $>_s$.
 We call $>_s$ the starting order for the walk, we will assume that we have some matrix $M_s$ with the first row
 $\mathbf{w}_{s}$ representing $>_s$. We have $S$ corresponds to a cone $C_{S,>_s}$ in the Sagbi fan of $\mathcal{A}$.
 The goal is to compute a Sagbi basis of $\mathcal{A}$ with respect to some target order $>_{t}$.
 This monomial order can be represented by by some matrix $M_t$ with first row $\mathbf{w}_{t}$.
 As we have both $>_s$ and $>_t$ are global orderings so $\mathbf{w}_s$ and $\mathbf{w}_t$ lie in positive orthant
 which is convex so we can use a straight line  between the two points, $(1-u)\mathbf{w}_{s}+u\mathbf{w}_{t}$ for
 $u\in [0,1]$. The Sagbi walk consist of two basic steps:\\ \\
$\bullet$ Crossing from one cone to next.\\
$\bullet$ Computing the Sagbi basis if $\mathcal{A}$ corresponding to the new cone.
\subsection{Crossing Cones}
In this section we discuss the procedure to cross from one cone to another of the Sagbi Fan, which is already given for
 the Gr\"{o}bner fan in \cite{CLO1}, Chapter 8. We can used it for the Sagbi fan because geometrically they are same. We overview  this procedure for Sagbi cones to used in Sagbi walk algorithm. For details see \cite{CLO1}, chapter 8, page 437.

Assume we have Sagbi basis $S_{old}$ corresponding to the cone $C_{old,>_{old}}$, and matrix $M_{old}$ with first row $\mathbf{w}_{old}$ representing $>_{s}$. As we continue along the path from $\mathbf{w}_{old}$, let $\mathbf{w}_{new}$ be the $last\, point$ on the path that lies in the cone $C_{old,>_{old}}$.

Let $S_{old}=\{\,x^{\alpha(i)}+\sum_{i,\beta}c_{i,\beta}x^{\beta}\,\,:\,\,1\leq i\leq t\,\}$, where $x^{\alpha(i)}$ is the leading term of $g_i$ with respect to $>_{M_{old}}$. To simplify notation, let $v_1,\ldots,v_m$ denote the vectors $\alpha(i)-\beta$ where $1\leq i\leq t$ and $c_{i,\beta}\neq 0$.
The new weight vector $\mathbf{w}_{new}$ is given by
 $\mathbf{w}_{new}=(1-u_{last})\mathbf{w}_{old}+u_{last}\mathbf{w}_{t}$ where $u_{last}$ is computed by the algorithm given in \cite{CLO1}, page 437. For the reader convenience we give this algorithm :
\begin{alg}
$\,\,$\\
Input : $\mathbf{w}_{old},\mathbf{w}_{t},v_1,\ldots,v_m$\\
Output : $u_{last}$
 \begin{itemize}
 \item $u_{last}=1$
 \item For $j=1,\ldots,m$ \\ \\
\indent If $(\mathbf{w}_{t}.v_{j}< 0)$;\\
\indent \quad Then $u_{j}=\frac{\mathbf{w}_{old}.vj}{\mathbf{w}_{old}.v_j-\mathbf{w}_t.v_j};$\\ \\
\indent  If $(u_{j}< u_{last})$\\
\indent \quad Then $u_{last=u_{j}}$;\\
\item return $u_{last}$;
\end{itemize}

\end{alg}

Once we have $\mathbf{w}_{new}$, we need to choose the next cone in the Sagbi fan. Let $>_{new}$ be the weight order
where we first compare $\mathbf{w}_{new}$-$degree$ and break ties using the target order. Since $>_t$ is represented by
$M_t$, it follows that $>_{new}$ is represented by $(\mathbf{w}_{new},M_{t})$. This gives the next cone $C_{new,>_{new}}$.
The following lemma shows that whenever $\mathbf{w}_{old}\neq \mathbf{w}_{t}$, the above process
guaranteed to move closer to $\mathbf{w}_{t}$.
\begin{lem}(c.f. \cite{CLO1}, chapter 8, page 438 )
Let $u_{last}$ be as in algorithm $1$ and assume that $>_{old}$ is represented by $(\mathbf{w}_{old},M_{t})$.
 Then $u_{last}> 0$.
\end{lem}
\subsection{Converting Sagbi Basis}
Once we have crossed from $C_{old}$ into $C_{new}$, we need to convert the the Sagbi basis $S_{old}$
(which are Sagbi basis of $\mathcal{A}$ with respect to $>_{old}$)  into Sagbi basis for $\mathcal{A}$
with respect to monomial order $>_{new}$, represented by $(\mathbf{w}_{new},M_{t})$. This can be done as follows.

The key feature of $\mathbf{w}_{new}$ that it is lies on the boundary of $C_{old,>_{old}}$, so that some
of the inequalities
become equalities. This mean that the leading term of some $g\in S_{old}$ has same $\mathbf{w}_{new}$-degree some other term
in $g$.
Note that $\mathbf{w}_{new}\in C_{old,>_{old}}$ guarantees that the $LT_{>_{old}}(g)$ appear in $In_{w}(g)$. Suppose that $H$ is
Sagbi basis of $in_{w_{new}}(S_{old})$ with respect to $>_{new}$. It is surprising  that once we have $H$, it is relatively
easy to convert $S_{old}$ into desired Sagbi basis.

\begin{thm}
Let $S_{old}=\{g_1,\ldots,g_{t}\}$ be the Sagbi basis for a subalgebra with respect to $>_{old}$.
 Let $>_{new}$ be represented by $(\mathbf{w}_{new},M_{t})$, where $\mathbf{w}_{new}$ is any weight vector in $C_{old}$ and
let $H=\{h_1,\ldots,h_s\}$ be the monic Sagbi basis of  $in_{w_{new}}(S_{old})$ with respect to $>_{new}$ as above. Express each $h_j\in H$ as
$$h_{j}=P_{j}(In_{w}(g_{1}),\ldots,In_{w}(g_{t})),\,\,\,\,P_{j}\in K[y_1,\ldots,y_t],\,\,g_{i}\in S_{old}  $$
Then replacing the initial terms by the $g_{i}$ themselves, the polynomials
$$\overline{h}_{j}=P_{j}(g_{1},\ldots,g_{t}),\,\,\,1\leq j\leq s $$
form a Sagbi basis of $\mathcal{A}$ with respect to $>_{new}$.
\end{thm}

 We say a weight vector $\mathbf{w}$ is compatible with a monomial $>$ if $LT_{>}(f) $ appears in $In_{w}(f)$ for all nonzero
polynomials $f$. Before proof of theorem 3.6 we will proof following lemma.
\begin{lem}
Fix $\mathbf{w}\in (\mathbb{R}^{n})^{+}\backslash \{0\}$ and let $S$ be the Sagbi basis of a subalgebra $\mathcal{A}$ for a monomial order $>$.\\
a. If $\mathbf{w}$ is compatible with $>$, then $LT_>(\mathcal{A})=LT_{>}(In_{w}(\mathcal{A}))=LT_{>}(in_{w}(\mathcal{A}))$.\\
b. If $\mathbf{w}\in C_{S,>}$, then $In_{w}(S)$ is Sagbi basis of $in_{w}(\mathcal{A})$ with respect to $>$. In particular,
 $$in_{w}(\mathcal{A})=in_{w}(S) $$
\end{lem}
\begin{proof}
For the part a, the first equality is $LT_>(\mathcal{A})=LT_{>}(In_{w}(\mathcal{A})) $ is obvious since the leading term of any $f\in K[x_1,\ldots,x_n]$ appear in $LT_{w}(f)$. For the second equality we first show that $LT_{>}(f)\in LT_>(In_{w}(\mathcal{A}))$ whenever $f\in in_{w}(\mathcal{A})$. Given such an $f$, write it as
$$f=q(In_{w}(f_{1}),\ldots,In_{w}(f_{t})),\,\,\,\,q\in K[y_1,\ldots,y_t],\,\,f_{i}\in \mathcal{A}$$
Each side is a sum of $w$-homogenous components. Since $In_{w}(f_{i})$ is already w-homogenous, so any $In_{w}(f_i)$-$monomial$ is also $w$-homogenous, this implies that
$$   In_{w}(f)=\overline{q}(In_{w}(f_{1}),\ldots,In_{w}(f_{t})), $$
where we can assume that all the terms in $\overline{q}(In_{w}(f_{1}),\ldots,In_{w}(f_{t}))$ are $w$-homogenous with same $w$-degree. It follows that $In_{w}(f)=In_{w}(\overline{q}(f_1,\ldots,f_t))\in In_{w}(\mathcal{A})$. Then compatibility implies that $LT_>(f)=LT_>(In_{w}(f))\in LT_{>}(In_{w}(\mathcal{A}))$. The other inclusion is obvious

Turning to part b, first assume that $w$ is compatible with $>$. Then
$$in_{>}(\mathcal{A})=in_{>}(S)=in_{>}(In_{w}(S)) $$
where the first equality follows since $S$ is Sagbi basis for $>$ and the second follows since $w$ is compatible with $>$. Combining with part a, we see that  $in_{>}(in_{w}(\mathcal{A}))=in_{>}(In_{w}(S))$ (since from part a we have $K[LT_>(\mathcal{A})]=K[LT_{>}(in_{w}(\mathcal{A}))]=in_{>}(in_{w}(\mathcal{A}))$). Hence $In_{w}(S)$ is Sagbi basis of $in_{w}(\mathcal{A})$ for $>$ and the final assertion of lemma follows.

It remain to show that what happens when $\mathbf{w}\in C_{S,>}$. Consider the ordering $>_{w}$,
which first compare the $\mathbf{w}$-degree and ties using $>$. Note that $\mathbf{w}$ is compatible with $>_w$.
The key observation is that since $\mathbf{w}\in C_{S,>}$, the leading term of each $g\in S$ with respect to $>_{w}$
 are same with $>$. Therefore it follows that $S$ is Sagbi
basis for $>_{w}$. Since $\mathbf{w}$ is compatible with $>_{w }$ earlier part of the argument implies that
$In_{w}(S)$ is as Sagbi basis of $in_{w}(\mathcal{A})$ for $>_{w}$. However, for each $g\in S$, $In_{w}(g)$
has same leading term with respect to $>$ and $>_{w}$.
Again we conclude that $In_{w}(S)$ is Sagbi basis of $in_{w}(\mathcal{A})$ for $>$.
\end{proof}
Now we give the proof of theorem  .

\begin{proof}
 We will give the proof in three steps. Since $>_{new}$ is represented by $(\mathbf{w}_{new},M_{t}$), $\mathbf{w}_{new}$ is compatible
 with $ >_{new}$. By the part a of lemma $(14)$, we obtain
 $$LT_{>_{new}}(\mathcal{A})=LT_{>_{new}}(in_{w_{new}}(\mathcal{A}))$$
 The second step is to observe that since $\mathbf{w}_{new}\in C_{old,>_{old}}$, the final assertion of part b of
 lemma (18)
 implies
 $$ in_{w_{new}}(\mathcal{A})=in_{w_{new}}(S_{old}) $$
 For the third step we show that
 $$ in_{w_{new}}(S_{old})=in_{>_{new}}(H)=in_{>_{new}}(\overline{H}),$$
 where $H=\{h_1,\ldots,h_{s}\}$ is the Sagbi basis of $in_{w_{new}}(S_{old})$ and
 $\overline{H}=\{\overline{h}_1,\ldots,\overline{h}_s\}$ as described in the statement of theorem. The first equality is obvious(as $w_{new}$ is compatible with $>_{new}$), and for the second it is suffices to show that for each j, $LT_{>_{new}}(h_j)=LT_{>_{new}}(\overline{h}_j).$ We have
  $$h_{j}=P_{j}(g_1,\ldots,g_t)=\sum_{\alpha^{(j)}}{(g_1)}^{{\alpha_{1}}^{(j)}}\ldots{(g_t)}^{{\alpha_{t}}^{(j)}}$$

 Observe that all the terms in ${(g_1)}^{{\alpha_{1}}^{(j)}}\ldots{(g_t)}^{{\alpha_{t}}^{(j)}}-{(In_{w}(g_1))}^{{\alpha_{1}}^{(j)}}\ldots{(In_{w}(g_t))}^{{\alpha_{t}}^{(j)}} $
 have smaller $w_{new}$-degree than  those  in the initial form of ${(g_1)}^{{\alpha_{1}}^{(j)}}\ldots{(g_t)}^{{\alpha_{t}}^{(j)}}$. Taking sum on all ${\alpha}^{j}$ we get
  $\widehat{h_j}=\overline{h}_{j}-h_{j}=P_{j}(g_{1},\ldots,g_{t})-P_{j}(In_{w}(g_{1}),\ldots,In_{w}(g_{t}))$.
  There we get $\overline{h}_j=h_j+\widehat{h_j}$ it shows that we get $\overline{h}_j$
  by adding terms of with smaller $w_{new}$-degree. Since $>_{new}$ is compatible with $w_{new}$, the added terms are smaller in the new order, so the leading term of $\overline{h}_j$ and $h_j$ with respect to $>_{new}$ is same.

  Combining the three steps, we obtain
  $$in_{>_{new}}(\mathcal{A})=in_{>_{new}}(\overline{H})$$
  Since $\overline{h_j}\in \mathcal{A}$ for all $j$, we conclude that $H$ is as Sagbi basis for $\mathcal{A}$ with respect to $>_{new}$ as claimed.
\end{proof}
\subsection{The Algorithm}
The following algorithm is a basic Sagbi walk, following the straight line segment from $\mathbf{w}_s$ to $\mathbf{w}_t$.
\begin{thm} Let\\ \\
$\mathbf{Next Cone}$ be the procedure that computes $u_{last}$ from algorithm $1$.
Recall that $\mathbf{w}_{new}=(1-u_{last})\mathbf{w}_{old}+u_{last}\mathbf{w}_{t}$ is the last weight vector along the path that lies
in the cone $C_{old_{>_{old}}}$ of the previous Sagbi basis $S_{old};$\\
$\mathbf{Lift}$ be the procedure that lifts a Sagbi basis for the $\mathbf{w}_{new}$-initial terms of the previous
Sagbi basis $S_{old}$ with respect to $>_{new}$ to the Sagbi basis $S_{new}$ following theorem $16$; and\\
$\mathbf{Interreduce}$ be the procedure that takes a given set of polynomials and interreduce
them with respect to the given monomial order.\\ \\
Then the following algorithm correctly computes a Sagbi basis for $\mathcal{A}$ with respect to $>_{t}$ and terminates in finitely steps:

\end{thm}
\begin{alg}
$\,\,$\\
\\Input: $M_{s}$ and $M_{t}$ representing the start and target order with first rows $\mathbf{w}_{s}$ and $\mathbf{w}_{t}$, $S_{s}$ = Sagbi basis with respect to $M_{s}$ .
\\Output:  $S_{new}$ = Sagbi basis with respect to $M_{t}$.
\begin{itemize}
\item $M_{old}:= M_{s} ;$
\item $S_{old}:= S_{s} ;$
\item $\mathbf{w}_{new}:= \mathbf{w}_{s} ;$
\item $M_{new}:=(\mathbf{w}_{new},M_{t})$ ;
\item $done:=false$ ;
\item while$(done=false)$ \\
\indent$In:=In_{\mathbf{w}_{new}}(S_{old}) $ ;\\
\indent $InS:=$Sagbibasis$(In,>_{M_{new}})$ ;\\
\indent $S_{new}:=$Lift$(InS,S_{old},In,M_{new},M_{old})$ ;\\
\indent $S_{new}:=$interreduce$(S_{new},M_{new})$ ;\\
\indent $u:=$NextCone$(S_{new},w_{new},\mathbf{w}_{t}) ;$\\
\indent if $(\mathbf{w}_{new}=\mathbf{w}_{t})$ then\\
\indent \quad $done=true$ ; \\
\indent else\\
\indent \quad  $M_{old}:=M_{new} ;$\\
\indent \quad $S_{old}:=S_{new}$ ;\\
\indent \quad $\mathbf{\mathbf{w}}_{new}:=(1-u)\mathbf{w}_{new}+u\mathbf{w}_{t}$ ;\\
\indent \quad $M_{new}:=(\mathbf{w}_{new},M_{t}) ;$
\item return $(S_{new})$ ;
\end{itemize}
\end{alg}
\begin{proof}
We traverse the line segment from $\mathbf{w}_{s}$ to $\mathbf{w}_{t}$. To prove the termination, observe that by theorem $13$,
the Sagbi fan of $\mathcal{A}=K[S_{s}]$ has only finitely many cones, each of which have finitely many bounding
hyperplanes as in definition $14$. Discarding these hyperplanes that contain the line segment from $\mathbf{w}_{s}$ to $\mathbf{w}_{t}$,
the remaining hyperplanes determine a finite set of distinguished points on our line segment.

Now consider $u_{last}$=Next cone$(S_{new},\mathbf{w}_{new},\mathbf{w}_{t})$ as the algorithm.
This uses Algorithm $1$ with $\mathbf{w}_{old}$ replaced by the current value of $\mathbf{w}_{new}$.
Further-more, notice that the monomial order always come from the matrix form $(\mathbf{w}_{s},M_{t})$.
It follows that hypothesis of lemma $17$ is always satisfied. If $u_{last}=1$, then the next value of $\mathbf{w}_{new}$
is $\mathbf{w}_{t}$, so the algorithm terminates after one more pass through the main loop.
On the other hand if $u_{last}=u_{j}< 1$, then the next value of $\mathbf{w}_{new}$ lies on the hyperplane
$\mathbf{w}.v_{j}=0$, which is one of our finitely hyperplanes.
However algorithm $1$ implies that $\mathbf{w}.v_{j}< 0$, and $\mathbf{w}_{new}.v_{j}\geq 0$,
so that hyperplane meets the line segment in single point. Hence the next value of $\mathbf{w}_{new}$
is one of  our distinguish point. Furthermore, lemma implies that $u_{last}> 0$ so that if the current
$\mathbf{w}_{new}$ differ from $\mathbf{w}_{t}$, then we must move to a distinguish point farther along the segment.
Hence we must eventually meet reach $w_{t}$, at which point algorithm terminates.

To prove the correctness, observe that in each pass through the main loop, the hypothesis of theorem $16$ are satisfied.
Furthermore once the value of $\mathbf{w}_{new}$ reaches $\mathbf{w}_{t}$, the next pass through main loop computes a Sagbi basis of
$\mathcal{A}$ for the monomial order represented by $(\mathbf{w}_{t},M_{t})$. It follow that the final value of $G_{new}$
is Sagbi basis for $>_{t}$.
\end{proof}

\begin{exmp}
Consider the subalgebra $\mathcal{A}=K[xy+z^{2},x^{2}y^{2}+y^{3}]$. For our convenience we underline the leading term of Sagbi bases with respect to a given monomial order.  We have $$S_{s}=\{z^{2}+xy,y^{3}+x^{2}y^{2}\} $$ is Sagbi basis of $\mathcal{A}$ with respect to lexicographical ordering induced by $z> y> x$. Suppose we want to determine the Sagbi basis with respect to lexicographical order induced by $x> y > z $. We could proceed as follows. Let\\

$M_{s}=$
\(
\left(
  \begin{array}{ccc}
    0 & 0 & 1 \\
    0 & 1 & 0 \\
    1 & 0 & 0 \\
  \end{array}
\right)
\)\\
so $\mathbf{w}_{s}=(0,0,1)$. similarly we have \\

$M_{t}=$
\(
\left(
  \begin{array}{ccc}
    1 & 0 & 0 \\
    0 & 1 & 0 \\
    0 & 0 & 1 \\
  \end{array}
\right)
\)\\
and $\mathbf{w}_{t}=(1,0,0)$. We will choose square matrices defining the appropriate monomial orders in all of the following computations by deleting appropriate linearly independent rows. We begin considering the order defined by \\
$M_{new}=$
\(
\left(
  \begin{array}{ccc}
    0 & 0 & 1 \\
    1 & 0 & 0 \\
    0 & 1 & 0 \\
  \end{array}
\right)
\)\\
(using the weight vector $\mathbf{w}_{new}=(0,0,1)$ first, then refining by the target order).
We have $In=\{z^{2},x^{2}y^{2}\}$. Sagbi basis of $K[In]$ with respect to $M_{new}$ is remain $\{z^{2},x^{2}y^{2}\}$.
Therefore   Sagbi basis of $A\mathcal{}$ with respect to $M_{new}$ does not change. We have
$$ S_{new}=\{z^{2}+xy,x^{2}y^{2}+y^{3}\} $$

We then call the next cone procedure(algorithm 1) with $\mathbf{w}_{new}$ in place of $\mathbf{w}_{old}$. The cone of $>_{M_{new}}$ is defined by the two inequalities obtained by comparing $z^{2}$ vs. $xy$ and $x^{2}y^{2}$ vs $y^{3}$. By algorithm 1
$u_{last}$ is the largest $u$ such that $(1-u)(0,0,1)+u(1,0,0) $ lies in the cone and is computed as follows:

$z^{2}$ vs. $xy$ :
$\hspace{1 cm} v_{1}=(-1,-1,2)$, $\mathbf{w}_{t}.v_{1}=-1< 0$ $\Rightarrow u_{1}=\frac{\mathbf{w}_{new}.v_{1}}{\mathbf{w}_{new}.v_{1}-(-1)}=\frac{2}{3}$\\

$  x^{2}y^{2} $ vs. $y^{3} $ :
$\hspace{1 cm} v_{2}=(2,-1,0)$, $\mathbf{w}_{t}.v_{2}=2\geq 0$ $\Rightarrow u_{2}=1 $\\ \\
Hence the new weight vector is $\mathbf{w}_{new}=(1-\frac{2}{3})(0,0,1)+\frac{2}{3}(1,0,0)=(\frac{2}{3},0,\frac{1}{3})$, and \\

$M_{new}=$
\(
\left(
  \begin{array}{ccc}
    \frac{2}{3} & 0 & \frac{1}{3} \\
    1 & 0 & 0 \\
    0 & 1 & 0 \\
  \end{array}
\right)
\)\\
are updated for the next pass through main loop.

In the second pass, $In=\{z^{2}+xy,x^{2}y^{2}\}$. We compute  the Sagbi basis for $K[In]$ with respect to $>_{new} $(with respect to this order, the leading term of the first element is $xy$), and find
$$H=\{h_1=xy+z^{2},h_2=x^{2}y^{2},h_3=xyz^{2}+\frac{1}{2}z^{4}\} $$
In terms of generators for $K[In]$, we have
$$xy+z^{2}=P_1(z^{2}+xy,x^{2}y^{2}),\,\,\,\, P_1\in K[t_1,t_2]\,\,\hbox{and}\,\,P_1(t_1,t_2,)=t_1  $$
$$ x^{2}y^{2}=P_2(z^{2}+xy,x^{2}y^{2}),\,\,\,\, P_2\in K[t_1,t_2]\,\,\hbox{and}\,\,P_2(t_1,t_2)=t_2   $$
$$xyz^{2}+\frac{1}{2}z^{4}=P_3(z^{2}+xy,x^{2}y^{2}),\,\,\,\, P_3\in K[t_1,t_2]\,\,\hbox{and}\,\,P_3(t_1,t_2)={t_1}^{2}-t_2 $$
So by theorem $16$, to get the next Sagbi basis we lift to
$$P_1(z^{2}+xy,x^{2}y^{2}+y^{3})=xy+z^{2}  $$
$$P_2(z^{2}+xy,x^{2}y^{2}+y^{3})=x^{2}y^{2}+y^{3}   $$
$$P_3(z^{2}+xy,x^{2}y^{2}+y^{3})= xyz^{2}+\frac{1}{2}z^{4}-\frac{1}{2}y^{3} $$
Interreducing with repect to $>_{new}$, we obtain the  Sagbi basis $S_{new}$ given by
$$\{ xy+z^{2}, xyz^{2}+\frac{1}{2}z^{4}-\frac{1}{2}y^{3}  \}$$
The call to NextCone returns $u_{last}=1 $, since there is no pair of terms equal weight for any point on the line segment  $(1-u)(\frac{2}{3},0,\frac{1}{3})+u(1,0,0) $. Thus $\mathbf{w}_{new}=w_{t}$, after one more pass through the main loop, during which $S_{new}$ doesnot change, the algorithm terminates. Hence the final output is
$$S_{t}= \{ xy+z^{2}, xyz^{2}+\frac{1}{2}z^{4}-\frac{1}{2}y^{3}  \}    $$
which is Sagbi basis of $\mathcal{A}$ with respect to target order.

\end{exmp}

 ${\mathbf{Acknowledgements}}.$ Thanks to my PhD supervisor Prof Dr Gerhard Pfister for his constant support and 
 valuable suggestions.


\end{document}